\documentclass[11pt]{article}
\usepackage{amssymb,amsmath,amsthm,mathrsfs,amsfonts,graphicx,epsf}
\usepackage{color}
\usepackage{xcolor}
\usepackage[active]{srcltx}

\usepackage{hyperref}            

 \usepackage{esint}

 \usepackage{fullpage}

\newcommand{\kdiff}{m-$\con^1_k$}

\newcommand{\meas}{\mathrm{meas}}
\newcommand{\bD}{{\cal D}}

\newcommand{\order}{degree}

\newcommand{\Tt}{\tilde}

\newtheorem{theorem}{Theorem}
\newtheorem{corollary}{Corollary}
\newtheorem{lemma}{Lemma}
\newtheorem{proposition}{Proposition}
\theoremstyle{definition}
\newtheorem{definition}{Definition}

\theoremstyle{remark}

\newtheorem{remark}{Remark}
\newtheorem{example}{Example}
\newcommand{\bt}{\begin{theorem}}
\newcommand{\et}{\end{theorem}}
\newcommand{\bl}{\begin{lemma}}
\newcommand{\el}{\end{lemma}}
\newcommand{\bp}{\begin{proposition}}
\newcommand{\ep}{\end{proposition}}
\newcommand{\bc}{\begin{corollary}}
\newcommand{\ec}{\end{corollary}}
\newcommand{\bdeff}{\begin{definition}}
\newcommand{\edeff}{\end{definition}}
\newcommand{\brem}{\begin{remark}}
\newcommand{\erem}{\end{remark}}
\renewcommand{\r}[1]{(\ref{#1})}
\newcommand{\con}{{\mathcal C}}

\def\cc{\mathit{C}}
\def\diam{\mathop{\mathrm{diam}}}
\def\ll{\mathcal{L}}
\def\ss{\mathcal{S}}
\newcommand{\hh}{{\mathcal H}}

\newcommand{\Len}{\mathrm{Length}}

\newcommand{\hd}{\widehat{d}}

\newcommand{\bi}{\begin{itemize}}
\newcommand{\iii}{\item}
\newcommand{\ei}{\end{itemize}}
\newcommand{\bd}{\begin{description}}
\newcommand{\ed}{\end{description}}

\newcommand{\bqn}{\begin{eqnarray}}
\newcommand{\eqn}{\end{eqnarray}}
\newcommand{\eqnn}{\nonumber\end{eqnarray}}

\newcommand{\ba}[1]{\begin{array}{#1}}
\newcommand{\ea}{\end{array}}

\newcommand{\R}{\mathbb{R}}

\newcommand{\N}{\mathbb{N}}

\newcommand{\g}{\gamma}
\newcommand{\al}{\alpha}
\newcommand{\eps}{\epsilon}

\title{\LARGE \bf
A new class of $(\hh^k,1)$-rectifiable subsets of metric spaces\thanks{This work was supported by the Digiteo grant {\it Congeo} and by the ANR project {\it GCM}, program ``Blanche'',
project number NT09\_504490.}
}
\author{
 R.~Ghezzi\thanks{Department of Mathematical Sciences and CCIB, Rutgers University 311 N 5$^\mathrm{th}$ Street  Camden, NJ 08102;
CMAP, \'Ecole Polytechnique
Route de Saclay, 91128 Palaiseau Cedex, France   and Team  GECO, INRIA Saclay --
\^Ile-de-France,    {\tt roberta.ghezzi@rutgers.edu} }
, F.~Jean\thanks{ENSTA ParisTech, UMA, 32, bd Victor, 75015 Paris, France  and Team  GECO, INRIA Saclay --
\^Ile-de-France,  {\tt frederic.jean@ensta-paristech.fr}}}

\begin{document}

\maketitle

\begin{abstract}
The main motivation of this paper arises from the study of Carnot--Carath\'eodory   spaces, where  the class of $1$-rectifiable  sets does not contain smooth non-horizontal curves; therefore a new definition of rectifiable sets including  non-horizontal curves is needed.
This is why we introduce  in any metric space a new class of curves,
 called continuously metric differentiable of \order\ $k$, which are H\"older but not Lipschitz continuous  when $k>1$. Replacing Lipschitz curves by this kind of curves we define $(\hh^k,1)$-rectifiable sets and show a density result generalizing the corresponding one in Euclidean geometry. This theorem is a consequence of computations of Hausdorff measures along curves, for which we give an integral formula.
 In particular, we show that both spherical and usual Hausdorff measures along   curves coincide with a class of dimensioned lengths and are related to an interpolation complexity, for which estimates have already been obtained in Carnot--Carath\'eodory   spaces.
\end{abstract}

\tableofcontents

\section{Introduction}

The main motivation of this paper arises from the study of Carnot--Carath\'eodory   spaces. Recall that such a metric space $(M,d)$ is defined by a sub-Riemannian manifold
 $(M,\bD,g)$, where  $M$ is a
smooth manifold, $\bD$ a subbundle of
$TM$ and $g$ a   Riemannian metric on $\bD$.
 The absolutely continuous paths which are
almost everywhere tangent to $\bD$ are called horizontal and their  length  is
obtained as in Riemannian geometry integrating the norm of
their tangent vectors.  The
distance $d$ is defined as the infimum of length of horizontal paths between two given points.

By construction, only horizontal paths may have finite length and may be   Lipschitz with respect to the distance. In contrast to the Euclidean case,
 both properties are independent on the regularity: all smooth non-horizontal paths have infinite length and are not Lipschitz. This gives rise to two kind of questions.

 The first query concerns the measure of non-horizontal paths: what kind of notion is the best suited? One of our motivation is that, from an intrinsic point of view, computing measures of paths should allow to determine metric invariants of curves and thus recover metrically the structure of the manifold~\cite{gro96}.
Since non-horizontal paths have a  metric dimension  greater than one (see~\cite{cocv_fj}), Hausdorff measures are the most natural candidates. However they pose  two problems: first they can hardly be computed (except for specific cases \cite{balu}), second they do not appear as integrals along the path, which  is what we expect for a measure  generalizing the notion of length.

The second question comes from geometric measure theory. A typical problem in this field is whether it is possible to characterize the geometric structure of a set using only measures. This gave rise to  the notion of rectifiable sets,	 
which is based on Lipschitz functions, and to density results in Euclidean   (see \cite{bes1,fed, preiss} and \cite{mattila} for a complete presentation) and general metric spaces  \cite{kircheim}. In the context of Carnot--Carath\'eodory spaces rectifiable sets have been studied in Heisenberg groups (see \cite{rect-heis,rect-heis2}) and a different notion of rectifiability was proposed in \cite{magnani}.
However,  in these spaces  the class of Lipschitz paths is quite poor and does not include  non-horizontal smooth curves which consequently are not rectifiable in the usual sense. To take into account the latter  curves we need  to define rectifiability through a larger class of paths, intrinsically characterized  by the distance.

\smallskip

In this paper we address these issues in any metric space, not only in Carnot--Carath\'eodory ones, by defining a class of curves in the spirit of (\cite{ambrosio-rectifiable,kircheim}). Namely, we introduce curves on a metric space $(M,d)$ that are continuously metric differentiable of \order\ $k$ (\kdiff\ for short) as  continuous curves $\g:[a,b]\to M$ such that the map
$$
t\mapsto \meas^k_t(\g)=\left(\lim_{s\to 0}\frac{d(\g(t+s),\g(t))}{|s|^{1/k}}\right)^k
$$
is well-defined and  continuous (see Definition~\ref{def:mdiff}). In an Euclidean space, this definition is useless since the class of \kdiff\ curves with non-zero measure is empty when $k>1$ (see Proposition~\ref{prop:eucl}). However, in the sub-Riemannian context, for integer values of $k$ this class of curves contains some smooth non-horizontal paths (see Proposition~\ref{ex:uno}).

For \kdiff\ curves  we can compute different kind of measures.
First, we examine the
Hausdorff measures: the usual ones $\hh^k$ and the spherical ones
$\ss^k$.  Second, we study the  $k$-dimensional length of a curve $\g:[a,b]\to M$ introduced
in~\cite{j-falbel}  and defined by
$$
\Len_k(\g([a,b]))=\int_a^b\meas^k_t(\g)dt.
$$
Third, we consider a measure based on approximations by
finite sets called interpolation complexity
(see~\cite{gau06,jea00}).
The first result of the paper (Theorem~\ref{th:main}) states that for an injective \kdiff\ curve $\g:[a,b]\to M$ we have
$$
\hh^k(\cc)=\ss^k(\cc)=\Len_k(\cc),
$$
 where $\cc=\g([a,b])$. It also provides a relation between
 $\hh^k(\cc)$  and the   interpolation complexity. On the one hand,
 Theorem~\ref{th:main} gives an integral formula for the Hausdorff
 measure.  On the other hand, it essentially implies that the
 considered measures are equivalent. 
%
%
Another interesting  property of injective \kdiff\ curves  with
non-zero $k$-dimensional measure is that the   $k$-dimensional density
of $\hh^k\lfloor_\cc$ exists and is constant along the curve (see
Proposition~\ref{p:main}). 

\smallskip

We  define  $(\hh^k,1)$-rectifiable sets as sets that are covered, up
to $\hh^k$-null sets, by  countable unions of \kdiff\ curves  (see
Definition~\ref{def:rect}). This notion  is modeled on the definition of
$(\hh^k,k)$-rectifiable sets in $\R^n$, which are sets
that are covered, up to $\hh^k$-null sets, by countable unions of
image of $C^1$ maps from $\R^k$ to $\R^n$.  Thanks to the properties of \kdiff\ curves, we show a density result for sets that are rectifiable according to our definition. Namely, the second main theorem of the paper (Theorem~\ref{th:fico}) states that if a set $S$ is $\hh^k$-measurable and satisfies $\hh^k(S)<+\infty$, then being $(\hh^k,1)$-rectifiable implies that the upper and lower densities of $\hh^k\lfloor_S$ are bounded by positive constants.

 Theorem~\ref{th:fico} is inspired by the  result proved in  Federer \cite[Th. 3.2.19]{federer}, which states that for a $\hh^k$ measurable subset  $E$ of the Euclidean $n$-space $(\hh^k,k)$-rectifiability implies that the measure $\hh^k\lfloor_E$ has  $k$-dimensional density equal to $1$ at $\hh^k$-almost every point of $E$.  The converse of this fact was proved for $k=1$ and for a general measure $\mu$ in  \cite{moore}. Much later, Preiss showed not only that the converse holds true for any $k$, but also a stronger result:  there exists a constant $c>1$ (depending only on $n$ and $k$)  such that if
 $$
 0<\limsup_{
 r \rightarrow 0^+
}\frac{\mu(E\cap B(x,r))}{r^k}\leq c \liminf_{
 r \rightarrow 0^+
}\frac{\mu(E\cap B(x,r))}{r^k}<+\infty, ~~\textrm{ for a.e. } x\in E,
 $$
 then $E$ is $(\mu, k)$-rectifiable.
 Our Theorem~\ref{th:fico} implies that an estimate of the type above is satisfied by $(\hh^k,1)$-rectifiable sets.  An open question is whether  an analogous of Preiss' result  still holds in non-Euclidean metric spaces with our definition of $(\hh^k,1)$-rectifiability.

 Another open problem is to show a Marstrand's type result (see
 \cite[Th.~1]{marstrand}) for $(\hh^k,1)$-rectifiable subsets, at
 least in Carnot--Carath\'eodory spaces. In Section~\ref{sec:examples}
 we construct \kdiff\ curves in sub-Riemannian manifolds having
 nonzero $k$-dimensional measure for integer values of $k\geq1$. When
 the curve is absolutely continuous, it is easy to see that being
 \kdiff\ with non-vanishing $k$-dimensional measure implies that $k$
 is an integer (see Corollary~\ref{cor:mas}). The question is whether
 such result holds true without assuming   absolute continuity. 

\medskip

The  structure of the paper is the following. In
Section~\ref{sec:mder} we give the definition of \kdiff\ curves in
metric spaces and construct them in Carnot--Carath\'eodory spaces. We
then study measures along curves. In Section~\ref{sec:nm} we recall
different notions of measures.  In Section~\ref{sec:ar} we show an
auxiliary result for \kdiff\ injective curves with nonzero
$k$-dimensional measure. In Section~\ref{sec:riem} we analyse \kdiff\
curves in  (the Euclidean space or a) Riemannian manifold. The main
theorem concerning injective \kdiff\ curves is proved in
Section~\ref{sec:mt}. 
%
%
Some possible generalizations  to non \kdiff\ curves are discussed
in Section~\ref{sec:tg}. 
Finally, in Section~\ref{sec:rect} we define
$(\hh^k,1)$-rectifiable sets and prove the density result.

\section{\kdiff\ curves}\label{sec:mder}

 Throughout the paper $(M,d)$ denotes a metric space.

\subsection{Definitions}

Let $\g:[a,b]\rightarrow M$ be a continuous curve, where $a,b\in\R$, and let $k\geq 1$ be a real number.
\begin{definition}\label{def:mdiff}
   We say that $\g$ is \emph{m-differentiable of \order\  $k$ at $t\in[a,b]$} if  the limit
\begin{equation}\label{eq:mder}
\lim_{\substack{s\rightarrow 0\\ t+s\in[a,b]}}\frac{d(\g(t+s), \g(t))}{|s|^{1/k}}
\end{equation}
exists and is finite. In this case, we call this limit the \emph{metric derivative of \order\ $k$ of $\g$ at $t$} and  we define moreover the \emph{$k$-dimensional infinitesimal measure of $\g$ at $t$} as
$$
\meas_t^k(\g)=\left(\lim_{\substack{s\rightarrow 0\\ t+s\in[a,b]}}\frac{d(\g(t+s),\g(t))}{|s|^{1/k}}\right)^k.
$$
When $\g$ is not m-differentiable of \order\ $k$ at $t$  we set $\meas_t^k(\g)=+\infty$.
\end{definition}
For the case $k=1$, the notion of metric derivative  is classical, see~\cite[Def. 4.1.2]{ambrosio}.  The $k$-dimensional infinitesimal measures of curves were introduced in the context of sub-Riemannian geometry in \cite{j-falbel}.

Note that if $\g$ is m-differentiable of \order\ $k$ at $t$ then, for any $k'$,
$$
\lim_{\substack{s\rightarrow 0\\ t+s\in[a,b]}}\frac{d(\g(t+s),\g(t))}{|s|^{1/k'}}=\lim_{\substack{s\rightarrow 0\\ t+s\in[a,b]}}\frac{1}{|s|^{1/k'-1/k}}\frac{d(\g(t+s),\g(t))}{|s|^{1/k}}.
$$
Therefore, for any  $k'>k$, $\meas_t^{k'}(\g)=0$. If moreover $\meas_t^k(\g)>0$, then for any $k'<k$
$\meas_t^{k'}(\g)=+\infty$.

\bdeff\label{def:kdiff}
Given  $k\geq1$, we say that $\g$ is  \emph{differentiable of class \kdiff\ on $[a,b]$}
(\kdiff\ for short) if  for every $t\in [a, b]$ the curve  is m-differentiable of \order\ $k$ at $t$ and the map $t\mapsto \meas^k_t(\g)$ is continuous.
\edeff
Clearly,  $\g$ is \kdiff\ if and only if the limit in \r{eq:mder}  exists and depends continuously on $t$.

We shall see in the next section that when a smooth structure on $M$ exists, \kdiff\ curves need not be differentiable in the usual sense. The following lemma states that they are H\"older continuous of exponent $1/k$ as functions from an interval to the metric space $(M,d)$.

\begin{lemma}
\label{le:holder}
Let $\gamma : [a,b] \to M$ be \kdiff\ on $[a,b]$, $k\geq 1$.
For any $t$ and $t+s$ in
$[a,b]$,
\begin{equation}\label{eq:unifest}
d(\gamma(t),\gamma(t+s))=|s|^{1/k}(\meas_t^k(\g)^{1/k}+\epsilon_t(s)),
\end{equation}
where  $\epsilon_t(s)$ tends to
zero as $s$ tends to zero uniformly with respect to $t$.
%
\end{lemma}
 This is a direct consequence of the continuity of $t\mapsto\meas_t^k(\g)$ on the compact interval $[a,b]$.

\subsection{Construction   of \kdiff\ curves}\label{sec:examples}

In this section we consider a class of metric spaces which are also smooth manifolds and construct smooth \kdiff\ curves on them with non-vanishing metric derivative of \order\ $k$ for some integer values of $k$. The analysis of this class of spaces is the main motivation of this paper.

Let $(M,d)$ be a metric space defined by a sub-Riemannian manifold $(M,\bD,g)$, i.e., $M$ is a
smooth  manifold, $\bD$ a
 subbundle  of $TM$, $g$ a  Riemannian metric on $\bD$, and $d$ is the associated sub-Riemannian distance.
We assume that Chow's condition is satisfied: let $\bD^s$ denote the
$\R$-linear span of brackets of degree $< s$ of vector fields tangent to
 $\bD^1 = \bD$; then, at every $p\in M$, there exists an integer $r = r(p)$
  such that $\bD^{r(p)}(p) = T_pM$, that is,
  \begin{equation}\label{eq:filtr}
\{0\}\subset \bD^1(p)\subset\bD^2(p)\subset\cdots \subset \bD^{r(p)}(p)= T_pM.
\end{equation}

Let $A \subset M$.
A point $p\in A$ is said \emph{$A$-regular} if the sequence
of dimensions $n_i(q)=\dim\bD^i(q),$ $i=1,\dots r(q)$
 remains
constant for $q \in A$ near $p$, and \emph{$A$-singular}
otherwise. The set  $A$ is said
\emph{equiregular} if every point of $A$ is $A$-regular. A curve $\g:[a,b]\to M$ is \emph{equiregular} if $\g([a,b])$ is equiregular.

\bp\label{ex:uno}
Let $\g:[a,b]\rightarrow M$ be an equiregular curve of class $\con^1$
and $k\in \N$  such that  $\dot\g(t)\in \bD^k(\g(t))$ for every
$t\in[a,b]$. Then $\g$ 
is \kdiff\  on $[a,b]$.

 If moreover $\dot\g(t)\notin\bD^{k-1}(\g(t))$ for a given $t\in[a,b]$ then  $\meas^k_t(\g)\neq 0$.
\ep
The proof of this proposition  is based on the notions of nilpotent approximation and privileged coordinates (see \cite{bellaiche}) and some results in \cite{j-falbel}.  We do not give  the complete argument, but only the underlying ideas. All the facts that here  are simply claimed are already established   and complete proofs can be found in the cited literature.

\smallskip

\noindent{\bf Sketch of the proof.} Since $\g([a,b])$ is equiregular, the integers $w_i$ defined by
$$
w_i = j, ~~\mathrm{ if }~ n_{j-1}(\g(t)) < i
\leq n_j(\g(t)),~~ i=1,\dots , n,
$$ do not depend on $t$. We  define for $s\geq 0$ the dilation $\delta_s:\R^n\rightarrow \R^n$ by
 $$
 \delta_s z = (s^{w_1} z_1,\dots , s^{w_n} z_n).
 $$
Moreover, locally there exist $n$
vector fields $Y_1, \dots, Y_n$ whose values at each $\g(t)$
form a basis of $T_{\g(t)}M$ adapted to the filtration \r{eq:filtr} at $\g(t)$,
in the sense that, for every integer $i \geq 1$, $Y_1(\g(t)),
\dots, Y_{n_i}(\g(t))$ is a basis of $\bD^i(\g(t))$. 
 The local diffeomorphism
$$
x \in \R^n \mapsto \exp (x_nY_n) \circ \cdots \circ \exp (x_1Y_1)
(\g(t))
$$
defines a system of coordinates $\phi^t: q  \mapsto
x=(x_1,\dots,x_n)$ on a neighborhood of $\g(t)$, satisfying $\phi^t(\g(t))=0$.  Following~\cite[Sec.~5.3]{bellaiche},  there exists a   sub-Riemannian distance
 $\hd_{t}$ on $\R^n$ such that
\begin{itemize}
\item $\hd_t$ is homogeneous under the dilation $\delta_s$, i.e., $\hd_t(\delta_s x,\delta_s x')= s \hd_t(x,x')$ for all $s\geq 0, x,x'\in\R^n$;
\item when defined, the mapping $t \mapsto \hd_t (\phi^t
  (q),\phi^t(q'))$ is continuous;
  \iii for $q$ in a neighborhood of $\g(t)$,  $
d(\g(t),q)=\hd_t(0,\phi^t(q))(1 +\eps_t(\hd_t(0,\phi^t(q))))
$, where  $\epsilon_t(s)$ tends to
zero as $s$ tends to zero uniformly with respect to $t$.
\ei
The coordinates  $\phi^t$  are privileged at $\g(t)$ and the distance
$\hd_t$ is the sub-Riemannian distance associated with a nilpotent
approximation at $\g(t)$.

Set $\phi^t(\g(t))=(\g_1(t),\dots,\g_n(t))$.
By the construction in the proof of  \cite[Le. 12]{j-falbel}, the limit
$$
\lim_{\substack{s\rightarrow 0\\
    t+s\in[a,b]}}\delta_{|s|^{-{1/k}}}\phi^t(\g(t+s))
$$
exists at every $t$ and is equal to
 $x(t)=(x_1(t),\dots,x_n(t))$, where
$$
x_j(t)=\left\{
\ba{ll}
0,& w_j\neq k\\
\dot \g_j (t),& w_j=k.
\ea
\right.
$$
Using the properties of $\hd_t$, we have
$$
\lim_{s\to0}\frac{d(\g(t+s),\g(t))}{|s|^{1/k}}=\lim_{s\to0}\frac{\hd_t(\phi^t(\g(t+s)),0)}{|s|^{1/k}}=\lim_{s\to 0}\hd_t(\delta_{|s|^{-{1/k}}}\phi^t(\g(t+s)), 0)=\hd_t(x(t),0).
$$
As a consequence, $\meas^k_t(\g)$ exists and is equal to $\hd_t(x(t),0)^k$.
Since the components of $x(t)$ are continuous and the distance $\hd_t$ depends continuously on $t$, $\g$ is \kdiff. If moreover $\dot\g(t)\notin\bD^{k-1}(\g(t))$ for a given $t\in[a,b]$ then $x(t)\neq 0$, whence $\meas^k_t(\g)\neq 0$.
\hfill$\blacksquare$

Let us explain the construction in Proposition~\ref{ex:uno} through an
example.  
\begin{example}
Consider the Heisenberg  group, that is, the sub-Riemannian manifold
$(\R^3,\bD, g)$ where $\bD$ 
 is the
linear span of the vector fields
$$
X_1(x,y,z)=(1,0,-y/2),~~X_2(x,y,z)=(0,1,x/2),
$$
 and $g=dx^2+dy^2$. Denote by $d$ the Carnot--Carath\'eodory distance
 associated with the Heisenberg group. Recall that $d$ is homogeneous
 with respect to the dilation  
 $$
 \delta_\lambda(x,y,z)=(\lambda x, \lambda y,\lambda^2 z),~~\lambda\geq 0.
 $$
and it is invariant with respect to the group law
$$
(x,y,z)*(x',y',z')=\left(x+x',y+y',z+z'+\frac12(xy'-x'y)\right).
$$ Moreover, for each point $(x,y,z)\in\R^3$, $\bD^2(x,y,z)=\R^3$ as
$[X_1,X_2](x,y,z)=(0,0,1)$.

 Let $\g(t)=(0,0,t)$, for $t\in\R$. Then $\g$ is of class m-$\con^1_2$
 and $\meas^2_t(\g)$ is a positive constant. 
 This is a consequence of Proposition~\ref{ex:uno} as $\g$ is smooth
 and, for all $t\in\R$,  $\dot\g(t)\in\bD^2(\g(t))$. 
 Let us compute explicitly $\meas^2_t(\g)$. Notice first that
 $d(\g(t+s),\g(t))=d((0,0,s),(0,0,0))$,  
    since $d$ is invariant with respect to the group law.
    Hence, using the homogeneity of $d$
    and the fact $d((0,0,1),0)=d((0,0,-1),0)$, 
$$
\lim_{s\to 0}\frac{d(\g(t+s),\g(t))}{|s|^{1/2}}=\lim_{s\to
  0}\frac{d((0,0,s),0)}{|s|^{1/2}}=d((0,0,1),0)= 2 \sqrt{\pi}, 
$$ 
the value of the distance resulting from an isoperimetric
problem. Note that such a computation can be generalized to any
contact 
sub-Riemannian manifold, see~\cite[Th.\ 22]{j-falbel}.
  \end{example}

\brem
Note that the equiregularity assumption is essential to obtain the
continuity  of $\hd_t$ and $\phi^t$ with respect to $t$. In particular
the proof of \cite[Le. 12]{j-falbel} is not valid without this
hypothesis\footnote{The statement of Lemma 12 in \cite{j-falbel} is
  incorrect. Indeed without equiregularity formula (3) therein does
  not hold.}. 
This assumption has also an intrinsic meaning. Indeed it is shown in
\cite{j-falbel} that the $k$-dimensional measure $\meas^k_t(\g)$ can
actually be defined through the distance  on the metric tangent space
to $(M,d)$ at $\g(t)$. Since the metric tangent space does not vary
continuously with respect to $t$ around $\cc$-singular points, where
$\cc=\gamma([a,b])$, in
general non equiregular curves may not be \kdiff.
\erem

Note that for every integer $k\in\{1,\dots, r(p)\}$, where $p$ is
regular, there exist $\con^1$ equiregular curves with tangent vector
belonging to $\bD^k\setminus\bD^{k-1}$.  As a consequence, for such
integers $k$ the class of \kdiff\ curves with non-vanishing metric
derivative of \order\ $k$ is not empty. For instance, this is the case
in the Heisenberg group for $k=2$, and in the Engel group (see below)
for $k=2,3$. On the contrary,  the next proposition states that  in
the Riemannian case, i.e., when $\bD=TM$, the class of \kdiff\ curves
with non-vanishing derivative is empty except for  $k=1$ (the proof of
Proposition~\ref{prop:eucl} is postponed to Section~\ref{sec:riem}).

\bp\label{prop:eucl}
Let $(M,g)$ be a Riemannian manifold. Let $k\geq 1$ and assume that
 $\g:[a,b]\to M$ is a \kdiff\ curve   such that
 $t\mapsto\meas^k_{t}(\g)$ does not vanish identically.
Then  $k=1$.
\ep

 Let $\g:[a,b]\to M$ be of class m-$\con^1_1$ and such that
 $\meas^1_t(\g)\neq 0$ for every $t\in[a,b]$. Then $\g$ is horizontal,
 i.e., it is absolutely continuous and $\dot \g(t)\in \bD^1(\g(t))$
 almost everywhere on $[a,b]$. To see this, remark that by
 construction,  $\g$ is  Lipschitz with respect to the sub-Riemannian
 distance.   The metric $g$ defined on $\bD$ can be extended  (at
 least in a tubular
 neighbourhood of $\gamma([a,b])$) to a
 Riemannian metric $\tilde g$ on $TM$. In this way we obtain a
 Riemannian distance on $M$ which is not greater than the
 sub-Riemannian distance. Hence $\g$ is  Lipschitz with respect to the
 chosen Riemannian distance which in turn implies that $\g$ is
 absolutely continuous. Therefore, by \cite[Pr. 5]{lipeq} $\g$ is
 horizontal, i.e., $\dot \g(t)\in \bD^1(\g(t))$ almost everywhere on
 $[a,b]$.

Using Proposition~\ref{ex:uno}, this fact can be partially generalized to
the case $k>1$ under the following form.
\bc\label{cor:mas}
Let $k\geq1$ and let $\g:[a,b]\to M$ be equiregular and of class
\kdiff, with $\meas^k_t(\g)\not\equiv 0$. If $\g$ is absolutely
continuous, then $k$ is the smallest integer $m$ such that $\dot \g(t)
\in \bD^m(\g(t))$ almost everywhere.
\ec
In particular Corollary~\ref{cor:mas} states that if
$\meas^k_t(\g)\not\equiv 0$ then $k$ is an integer, provided that $\g$
is absolutely continuous.  An open question is whether the latter
condition is necessary. If this were not the case then we would obtain
a Marstrand's type Theorem \cite[Th. 1]{marstrand} for \kdiff\ curves:
indeed we shall see in Proposition~\ref{p:main} that along injective
\kdiff\ curves with non-vanishing $k$-dimensional measure  the density
of $\hh^k$   exists and is constant.

Nevertheless, a \kdiff\  curve need not be $\con^1$ in the usual sense as it is shown below.

 \begin{example}\label{ex:due} Consider the Engel  group, that is, the sub-Riemannian manifold $(\R^4,\bD, g)$ where $\bD$ is the
linear span of the vector fields
$$
X_1(x,y,z,w)=(1,0,0,0),~~X_2(x,y,z,w)=(0,1,x, x^2/2),
$$
and $g=dx^2+dy^2$.
 Let
$\g(t)=(0,0,W(t), \varphi(t)),$ where  $\varphi\in \con^1$ and $W$ is the Weierstrass function $$
W(t)=\sum_{n=0}^{\infty}\al^n(\cos(\beta^n \pi t)-1),~~t\in\R,
$$
where $0 < \al < 1,$  $\beta > 1$, and $\al \beta >1$ see \cite{weier}. It was proved in \cite{hardy} that
    $W(t)$ is continuous,  nowhere
differentiable on the real line, and satisfies
\begin{equation}\label{eq:holder}
W(t+h)-W(t)=O(|h|^\xi),~~\xi=\frac{\log (1/\al)}{\log \beta}<1,
\end{equation}
uniformly with respect to $t \in\R$.
 Then, choosing $\al, \beta$ such that $\xi >2/3$, $\g$ is
  continuous and m-$\con^1_3$, but nowhere differentiable.
  Indeed, it is not hard to verify that the sub-Riemannian distance $d$ satisfies the following homogeneity property
  $$
\lambda d((0,0,\bar z, \bar w),(0,0,z,w))=d(0,(0,0,\lambda^2(z-\bar z),\lambda^3(w-\bar w))),
  $$
 for every $\lambda \geq 0$. Then
  we have
$$
\lim_{s\to 0}\frac{1}{|s|^{1/3}}d((0,0,W(t),\varphi(t)),(0,0,W(t+s),\varphi(t+s)))=\lim_{s\to 0}d\left(0,
\left(0,0,\frac{O(|s|^{\xi})}{|s|^{2/3}}, \varphi'(t)\right)\right). 
$$
Since $\xi>2/3$,  $\g$ is m-differentiable of \order\ $3$ at each $t$ and  $\meas^3_t(\g)=d(0,(0,0,0,\varphi'(t))^3$, which is non-zero for  a suitable choice of $\varphi$. Therefore,  $\g$ is m-$\con^1_3$  and by the properties of $W(t)$, $\g$ is nowhere differentiable.
\end{example}
Notice that if $\g$ is \kdiff\ and $k'\geq k$, then $\g$ is m-$\con^1_{k'}$. Define
$
k_\g\geq 1
$
as the infimum of $k\geq 1$ such that $\g$ is \kdiff.
Then $k_\g$ need not be an integer as it is shown in the next example. Moreover, $\g$ is not necessarily m-$\con^1_{k_\g}$.

\begin{example}
Consider the sub-Riemannian structure of Example~\ref{ex:due} and
 the curve $\g(t)=(0,0,W(t),0)$. Then $k_\g=2/\xi$ may be  any real number greater than $2$ (see \r{eq:holder}), but  $\g$ is not m-$\con^1_{k_\g}$.
\end{example}

\section{Measures along curves}\label{sec:meas}
This section is devoted to compute Hausdorff (and spherical Hausdorff) measures of
 continuous curves and to establish a relation with the
 $k$-dimensional length and with the complexity.

\subsection{Different notions of measures}\label{sec:nm}

Denote by $\diam S$ the diameter of a set $S \subset M$.
Let $k \geq 0$ be a real number. For every set $A \subset M$, we
define the
\emph{$k$-dimensional Hausdorff measure} $\hh^k$ of $A$  as $\hh^k(A)
= \lim_{\eps \to 0^+} \hh^k_\eps(A)$, where
$$
\hh^k_\eps(A) = \inf \left\{ \sum_{i=1}^\infty  \left(\diam S_i\right)^k
\, : \, A \subset \bigcup_{i=1}^\infty S_i, \  \diam S_i \leq \eps,
\ S_i \hbox{ closed set} \right\},
$$
and the \emph{$k$-dimensional spherical Hausdorff measure} $\ss^k$  of
$A$  as $\ss^k(A)
= \lim_{\eps \to 0^+} \ss^k_\eps(A)$, where
$$
\ss^k_\eps(A) = \inf \left\{ \sum_{i=1}^\infty  \left(\diam
S_i\right)^k \, : \, A \subset \bigcup_{i=1}^\infty S_i, \ S_i \hbox{ is
  a ball}, \ \diam
S_i \leq \eps  \right\}.
$$

In the Euclidean space $\R^n$,  $k$-dimensional Hausdorff measures
are often
defined as $2^{-k}\alpha(k)\hh^k$ and $2^{-k}\alpha(k)\ss^k$, where
$\alpha(k)$ is defined from the
usual gamma function as $\alpha(k) = \Gamma (\frac{1}{2})^k / \Gamma
(\frac{k}{2}+1)$. This normalization
factor is necessary for the $n$-dimensional Hausdorff measure and
the Lebesgue measure  coincide on $\R^n$.

For a given set $A \subset M$, $\hh^k(A)$ is a decreasing function
of $k$, infinite when $k$ is smaller than a certain value, and zero
when $k$ is greater than this value. We call \emph{Hausdorff
dimension} of $A$ the real number
$$
\dim_\hh A = \sup \{ k \, : \, \hh^k (A) = \infty \} = \inf \{ k \,
: \, \hh^k (A) = 0 \}.
$$
Note that $\hh^k \leq \ss^k \leq 2^k \hh^k$, so the Hausdorff
dimension can be defined equally from Hausdorff or spherical
Hausdorff measures.

When the set $A$ is a curve, another kind of dimensioned measures can be obtained from the integration of $k$-dimensional infinitesimal measures.
Let $\gamma: [a,b] \to M$ be a  continuous curve and
$\cc=\gamma([a,b])$. For $k \geq 1$, we define the \emph{$k$-dimensional
length} of $\cc$ as
\begin{equation}
\label{eq:deflenk}
\Len_k ( \cc ) = \int_a^b  \meas_t^k (\gamma) \,  dt.
\end{equation}
where $\meas_t^k (\gamma)$ is as in Definition~\ref{def:mdiff} (these lengths were introduced in \cite{j-falbel} in the sub-Riemannian context).
Thanks to the properties of $\meas_t^k(\g)$, $\Len_k(\g)$ is a decreasing function of $k$, infinite when $k$ is smaller than a certain value, and zero
when $k$ is greater than this value. We call this value the {\it length dimension of $\cc$}.

Another way to measure the set $\cc$  is to study its approximations by
 finite sets (see~\cite{jea00} and~\cite[p.\ 278]{gro96}). Here we only consider approximations by $\eps$-chains of $\cc$, i.e., sets of points $q_1=\gamma (a)$,
\dots, $q_N=\gamma(b)$ in $\cc$ such that $d(q_i,q_{i+1})\leq \eps$.
The \emph{interpolation complexity} $\sigma_{\mathrm{int}} (\cc, \eps)$ is the
minimal number of points in an $\eps$-chain of $\cc$.
This complexity has been computed in several cases in \cite{gz-nuovo}.

\brem
Notice that   for any injective m-$\con^1_1$ curve the equality   $\hh^1(\cc)=\Len_1(\cc)$  holds (see \cite[Th. 4.1.6, 4.4.2]{ambrosio}).
\erem

\subsection{\kdiff\  curves with non-vanishing $k$-dimensional measure}\label{sec:ar}

In this section we prove the following proposition about   \kdiff\ curves with non-vanishing $k$-dimensional measure. This result is the first step to prove Theorem~\ref{th:main}.

\bp\label{p:main}
Let $\g:[a,b]\rightarrow M$ be an injective \kdiff\ curve  and
$\cc=\g([a,b])$. Assume $\meas^k_t(\g)\neq 0$ for every $t$. Then
\begin{eqnarray}
 &&{\cal H}^k(\cc) = {\cal S}^k(\cc)=  \mathrm{Length }_k
 (\cc)\label{eq:measure}\\ 
 &&\lim_{\eps\rightarrow
   0^+}\eps^k\sigma_{\mathrm{int}}(\cc,\eps)=\Len_k(\cc), \label{eq:comp=len} 
\end{eqnarray}
and for every $q\in\cc$
\begin{equation}\label{eq:den}
\lim_{r\to 0^+}\frac{\hh^k(\cc\cap B(q,r))}{2 r^k}=1.
\end{equation}
\ep

\brem\label{rk:gen}
 Equations~\r{eq:measure}, \r{eq:den}   hold when we replace $[a,b]$
 by the open interval $(a,b)$. Also, they hold for unbounded
 intervals. Therefore, thanks to the regularity  of ${\cal L}^1$ and
 $\hh^k$ measures, equations ~\r{eq:measure}, \r{eq:den} are still
 verified when we replace $\cc$ by $\g(A)$, for any measurable set
 $A\subset [a,b]$.
 \erem
If we drop the injectivity assumption we obtain the following weaker result.
\bc\label{cor:leq}
 Let $\g:[a,b]\rightarrow M$ be a  \kdiff\ curve  and
$\cc=\g([a,b])$. Assume $\meas^k_t(\g)\neq 0$ for every $t$. Then
\begin{equation}\label{eq:leq}
\hh^k(\cc)=\ss^k(\cc)\leq\Len_k(\cc).
\end{equation}
\ec

\noindent{\bf Proof of Corollary~\ref{cor:leq}.} Since
$\meas^k_t(\g)\neq 0$ for every $t\in [a,b]$, $\g$ is  locally
injective. Hence $[a,b]$ is the disjoint union of a finite family of
intervals $I_i$ such that $\g|_{I_i}$ is injective. For each $i$ there
exists  a measurable subset  $A_i\subset I_i$  such that
$\cc=\cup\g(A_i)$ and the sets $\g(A_i)$ are pairwise disjoint. Using
Remark~\ref{rk:gen}, formula \r{eq:measure} applies to each
$\g(A_i)$. Since $\hh^k(\cc)=\sum_i \hh^k(\g(A_i))$ and
$\ss^k(\cc)=\sum_i \ss^k(\g(A_i))$, we obtain \r{eq:leq}.

\hfill$\blacksquare$

\smallskip


The proof of Proposition~\ref{p:main} is based on the following  result for bi-H\"{o}lder continuous curves.

\bl
\label{le:bounds}
Let $\g:[0,T]\rightarrow M$ be an injective curve  and
$\cc=\g([0,T])$.
Assume that  there
exist positive constants $\delta_-, \delta_+,$ and $\eta$ such that
\begin{equation}
\label{eq:holder12}
\delta_-|s|^{1/k} \leq d \big( \gamma(t),\gamma (t+s)
\big) \leq
\delta_+|s|^{1/k},
\end{equation}
for every $t,t+s \in [0,T]$ with  $|s|<\eta$. Then
\begin{eqnarray}
& \displaystyle
\delta_-^kT \leq {\cal
  H}^k(\cc) \leq
\delta_+^k T, \label{eq:measure2}& \\
& \displaystyle \delta_-^k T \leq \liminf_{\eps \to
  0^+} \eps^k \sigma_{\mathrm{int}} (\cc, \eps) \leq
 \limsup_{\eps\rightarrow
   0^+}\eps^k\sigma_{\mathrm{int}}(\cc,\eps) \leq \delta_+^k
 T,& \label{eq:measure3} \\
& \displaystyle   \ss^k (\cc) \geq \hh^k (\cc) \geq \left(\frac{\delta_-}{\delta_+}\right)^{2k} \ss^k (\cc), & \label{eq:measure4}
\end{eqnarray}
and, for every $t\in [0,T]$ and $r>0$ small enough,
\begin{eqnarray}
& \displaystyle  \left(\frac{\delta_-}{\delta_+}\right)^{k} \leq \frac{\hh^k (\cc \cap B(\gamma(t)
,r))}{2 r^k} \leq \left(\frac{\delta_+}{\delta_-}\right)^{k}. & \label{eq:measure5}
\end{eqnarray}
\el

\noindent{\bf Proof.}
Let $\eps >0$ be smaller than $\delta_+\eta^{1/k}$. We denote by
$N$ the smallest integer such that $T \leq N
(\frac{\eps}{\delta_+})^k$ and define $t_0, \dots, t_N$ by
$$
t_i = i \big(\frac{\eps}{\delta_+}\big)^k \quad \hbox{for } i= 0,
\dots, N-1, \qquad t_N = T.
$$
Set $S_i=\gamma([t_{i-1},t_i])$, $i=1, \dots , N$. For $t,t'$ in
$S_i$, one has $|t-t'|\leq \eps^k /\delta_+^k$; it follows
from~(\ref{eq:holder12}) that
\begin{equation}
\label{eq:d<eps}
  d( \gamma(t),\gamma(t')) \leq \delta_+ |t-t'|^{1/k} \leq \eps,
\end{equation}
which in turn implies $\diam S_i \leq \eps$. Thus $\hh^k_\eps
(\cc) \leq \sum_i (\diam S_i)^k \leq N \eps^k$. Using $(N-1)\eps^k
< T \delta_+^k$, we obtain
$$
\hh^k_\eps (\cc) \leq \delta_+^k T + \eps^k.
$$
It also results from inequality~\eqref{eq:d<eps}  that $\gamma(t_0)$, \dots,
$\gamma(t_N)$ is an
$\eps$-chain of $\cc$  which implies
$$
\eps^k \sigma_{\mathrm{int}} (\cc, \eps) \leq \delta_+^k T + \eps^k,
$$
Taking the limit as $\eps \to 0$ in
the preceding  inequalities, we find
$$
\hh^k (\cc)\leq \delta_+^k  T\quad
\mathrm{ and } \quad \limsup_{\eps \to
  0^+} \eps^k \sigma_{\mathrm{int}} (\cc, \eps) \leq
\delta_+^k T.
$$

We now  prove converse inequalities for $\hh^k$ and
$\sigma_{\mathrm{int}}$. 
Fix $\eps>0$ and  consider a countable family $S_1, S_2, \dots$ of
closed subsets of $M$
such that $\cc \subset \bigcup_i S_i$ and $\diam S_i \leq \eps$.  For
every $i\in \mathbb{N}$, we set
$I_i = \gamma^{-1} (S_i \cap \cc)$. As $\g$ is injective, if $\eps$ is
small enough then it results from~(\ref{eq:holder12}) that for any
$t,t'$ in $I_i$ there holds
$$
\diam S_i \geq d(\gamma (t),\gamma(t')) \geq \delta_- |t -
t'|^{1/k},
$$
which implies $\ll^1 (I_i) \leq (\diam S_i)^k / \delta_-^k$.
Note that $T \leq \sum_i \ll^1
(I_i)$ since the sets $I_i$ cover $[0,T]$. It follows that $\hh^k_\eps
(\cc) \geq T \delta_-^k$, that is,
\begin{equation}
\label{eq:hkgeq2}
\hh^k_\eps (\cc) \geq \delta_-^k T.
\end{equation}
In the same way, an $\eps$-chain $\gamma( t_0)=\gamma(0)$, \dots,
$\gamma(t_N)=\gamma(T)$ of $\cc$ satisfies $N \eps^k \geq T
\delta_-^k$ since the injectivity of $\g$ assures that
$$
\eps \geq d (\gamma( t_{i-1}), \gamma ( t_i)) \geq \delta_- |t_i -
t_{i-1}|^{1/k}.
$$
It follows that $\eps^k \sigma_{\mathrm{int}} (\cc, \eps) \geq
\delta_-^k T$.
Taking the limit as $\eps \to 0$ in
this inequality and in~(\ref{eq:hkgeq2}), we find $\hh^k (\cc)
\geq \delta_-^k T$ and
$\liminf_{\eps \to
  0^+} \eps^k \sigma_{\mathrm{int}} (\cc, \eps) \geq
\delta_-^k T$, which completes the proof of~\eqref{eq:measure2} and~\eqref{eq:measure3}.
\medskip

The first inequality in~\eqref{eq:measure4} always holds. Before proving the second one, let us recall a standard result in geometric measure theory (see for
instance~\cite[2.10.18, (1)]{federer}).
\emph{Let $X$ be a metric space and  $\mu$ be a regular measure on $X$ such
 that the closed balls in $X$ are $\mu$-measurable. If
$$
\limsup_{
 \substack{y\in B(x,r)\\r \rightarrow 0^+}
}\frac{\mu(B(y,r))}{(\diam B(y,r))^k} \geq \lambda,
$$
for every point $x\in X$,
then $\mu(X) \geq \lambda \ss^k(X)$.
}
We will apply this result to the metric space $(\cc,d|_{\cc})$ and to
the measure 
$\mu=\hh^k\lfloor_{\cc}$.

If $t\in [0,T]$, then for  $r>0$ small enough, there holds
\begin{equation}
\begin{array}{c}
\label{eq:ball_estimate}
 \gamma([t - \frac{r^k}{\delta_+^k} ,t + \frac{r^k}{\delta_+^k}])
\subset \cc \cap B(\gamma(t) ,r) \subset \gamma([t -
\frac{r^k}{\delta_-^k} ,t + \frac{r^k}{\delta_-^k}]).
\end{array}
\end{equation}
The diameter of this set then satisfies
$$
\diam (\cc \cap B(\gamma(t) ,r)) \leq d\left( \g\left(t - \frac{r^k}{\delta_+^k}\right),\g\left(t + \frac{r^k}{\delta_+^k}\right)\right) \leq \delta_+ \frac{2^{\frac{1}{k}}r}{\delta_-}.
$$
Moreover, applying~\eqref{eq:measure2} to the curve $\gamma$ restricted to $[t - \frac{r^k}{\delta_+^k} ,t + \frac{r^k}{\delta_+^k}]$, we obtain
$$
\hh^k (\cc \cap B(\gamma(t) ,r))
\geq \delta_-^k \frac{2r^k}{\delta_+^k}.
$$
Thus we have, for every point $\g(t')\in \cc$,
$$
\limsup_{
 \substack{r \rightarrow 0^+\\ \g(t)\in B(\g(t'),r)}
}\frac{\hh^k\lfloor_{\cc} (B(\gamma(t) ,r))}{(\diam (\cc \cap B(\gamma(t) ,r)))^k} \geq \left(\frac{\delta_-}{\delta_+}\right)^{2k},
$$
which implies $\hh^k (\cc) \geq \left(\frac{\delta_-}{\delta_+}\right)^{2k} \ss^k (\cc)$.
\medskip

Finally, formula~\eqref{eq:measure5} results from~\eqref{eq:measure2} applied to the restrictions of $\gamma$ in~\eqref{eq:ball_estimate}.
 
\hfill$\blacksquare$



\noindent{\bf Proof of Proposition~\ref{p:main}.}
By definition,
$$
\Len_k (\cc) = \int_a^b \meas_t^k(\g) \, dt.
$$
Note that the $k$-dimensional length does not
depend on the parameterization~\cite[Le. 16]{j-falbel}. Thus, up to a
reparameterization by the
$k$-length, we assume that $\g$ is defined on the interval $[0,T]$,
with $T=\Len_k (\cc)$, and that $\meas_t^k(\g)\equiv 1$.

Fix $\delta >0$. Then, by
Lemma~\ref{le:holder}, there exists $\eta >0$  so that the hypothesis
of Lemma~\ref{le:bounds} is satisfied with $\delta_-=1-\delta$ and
$\delta_+=1+\delta$. We let $\delta$ tends to zero
in inequalities~\eqref{eq:measure2}--\eqref{eq:measure5} and the proposition follows.
\hfill$\blacksquare$

\brem
Another way to measure $\cc$ using approximations by
finite sets is to consider $\eps$-nets,   i.e.,  sets of points $q_1,\dots, q_n\in M$ such that the union of closed balls $B(q_i,\eps)$ covers $\cc$,  and the metric entropy $e(\cc,\eps)$ which  is the minimal number of points in an $\eps$-net of $\cc$. Under the assumptions of Proposition~\ref{p:main},  the following estimates can be deduced for a \kdiff\ curve:
$$
\frac{\ss^k(\cc)}{2^k}\leq\liminf_{\eps \to 0^+} \eps^k e(\cc,\eps)\leq \limsup_{\eps \to 0^+} \eps^k e(\cc,\eps)\leq \frac{\ss^k(\cc)}{2}.
$$
\erem

\subsection{The Riemannian  case}\label{sec:riem}

Let us come back to the case where $(M,d)$ is a Carnot--Carath\'eodory space associated with a sub-Riemannian manifold $(M,\bD,g)$.
A consequence of Proposition~\ref{p:main} is that if  the structure is Riemannian, i.e., $\bD=TM$,  then the class of \kdiff\ curves having non-zero metric derivative of \order\ $k$ is empty if $k>1$, as stated in Proposition~\ref{prop:eucl}.

\medskip

\noindent{\bf Proof of Proposition~\ref{prop:eucl}.}
 Let $s\in[a,b]$ such that $\meas^k_s(\g)\neq 0$. Thus, restricted to a small enough neighbourhood $I=[s-\delta,s+\delta]$ of $s$, the curve $\g$ is injective and $\meas^k_t(\g)\neq 0$ on $I$. Up to  reparameterizing $\gamma|_{I}$, we may assume moreover $\meas^k_t(\g)\equiv 1$ on $I$. Also, it is sufficient to consider the case $M=\R^n$ and $d$ is the Euclidean distance on $\R^n$. 
Denote by $\cc$ the set $\g(I)$. By Proposition \ref{p:main}, we have $0<\hh^k(\cc)<+\infty$ and for every $t\in I$
$$
\lim_{r\to 0}\frac{\hh^k(\cc\cap B(\g(t),r))}{2 r^k}=1.
$$
Moreover, by Lemma~\ref{le:holder}, there exist
$0<\rho < 1$ such that 
$$
(1-\rho)|t-t'|^{1/k}\leq ||\g(t)-\g(t')||\leq (1+\rho)|t-t'|^{1/k}, \quad \forall\, t, t'\in I.
$$
The proposition then results from the lemma below.
 \hfill$\blacksquare$

\begin{lemma}
\label{le:riem_den}
Let $k\geq 1$ and let
 $\g:[a,b]\to \R^n$ be a bi-H\"{o}lder curve of exponent $1/k$, i.e.,
\begin{equation}\label{eq:bilip}
\delta_-|t-t'|^{1/k}\leq ||\g(t)-\g(t')||\leq \delta_+|t-t'|^{1/k},
\quad \forall\, t, t'\in [a,b], 
\end{equation}
with $\delta_-,\delta_+\neq 0$. Set $\cc=\g([a,b])$. If $0<\hh^k(\cc)<+\infty$
and if there exits a positive constant $c$ such that for every $t\in [a,b]$ 
$$
\lim_{r\to 0}\frac{\hh^k(\cc\cap B(\g(t),r))}{ r^k}= c,
$$ 
then  $k=1$.
\end{lemma}

\noindent{\bf Proof.} Under the assumptions of the lemma,  Marstrand's
Theorem \cite[Th. 1]{marstrand} 
assures that 
$k\in\N$ and $k\in\{1,2,\dots, n\}$.
Applying Preiss' result \cite{preiss}, there exists a countable family
of $k$-dimensional submanifolds $N_i \subset \R^n$ such that
$\hh^k(\cc\setminus\cup_{i}N_i)=0$.  Since  $\hh^k(\cc)>0$, there
exists $i$ such that $\hh^k(\cc\cap N_i)>0$. Let us rename $N_i$ by
$N$. Then  ${\cal L}^k_{N}(\cc\cap N)=\hh^k\lfloor_{N}(\cc)>0$, where
${\cal L}^k_{N}$ is the $k$-dimensional Lebesgue measure on $N$.
Hence there exists a density point $\g(t_0)\in\cc\cap N$, that is, a
point such that
$$
\lim_{r\to 0}\frac{{\cal L}^k_N(\cc\cap N\cap B_N(\g(t_0),r))}{{\cal
    L}^k_N(B_N(\g(t_0),r))}=1,
$$
where $B_N(\g(t_0),r)$ is the open ball in $N$ centered at $\g(t_0)$
of radius $r$.

 Let us identify $N$ with $\{(x_1,\dots,x_n)\in \R^n :
 x_{k+1}=\dots=x_n=0\}$ by choosing local coordinates around
 $\g(t_0)$. Using the inequalities \r{eq:bilip} and the density point
 $t_0$ it is not hard to prove that there exists a bi-Lipschitz
 homeomorphism from $B_N(\g(t_0),\delta_-)$ endowed with the Euclidean
 distance to $(-1,1)$ endowed with the distance $|\cdot|^{1/k}$ (see for instance the argument in the proof of  \cite[Pr. 4.12]{assouad}). Since the topological dimension of $B_N(\g(t_0),\delta_-)$ is $k$, then $k$ must be equal to $1$.

 \hfill$\blacksquare$

\subsection{Comparison of measures for  \kdiff\ curves}\label{sec:mt}

Next theorem generalizes the first part of Proposition~\ref{p:main} to
the case when the metric derivative may vanish. Namely, 
it compares the $\hh^k$  measure and the $\ss^k$  measure of  sets
that are  images of \kdiff\ curves. Also, it provides a relation among
such measures, the $k$-length  and the behaviour of the complexity of
the curve.

\bt\label{th:main}
Let $\g:[a,b]\rightarrow M$ be an injective \kdiff\ curve  and
$\cc=\g([a,b])$. Then
$$
 {\cal H}^k(\cc) = {\cal S}^k(\cc)=  \mathrm{Length }_k (\cc)=
 \lim_{\eps\rightarrow 0}\eps^k\sigma_{\mathrm{int}}(\cc,\eps).
$$
If moreover $\hh^k(\cc)>0$ or $\Len_k(\cc)>0$, then  for every $k'\geq 1$
$$
 {\cal H}^{k'}(\cc) = {\cal S}^{k'}(\cc)=  \mathrm{Length }_{k'} (\cc)=
 \lim_{\eps\rightarrow 0}\eps^{k'}\sigma_{\mathrm{int}}(\cc,\eps).
$$
\et
\brem
When $M$ is a sub-Riemannian manifold, in many cases Gauthier and
coauthors  (see \cite{gz-nuovo} and references therein) computed the
interpolation complexity of curves as integral of some geometric
invariants. Jointly with Theorem~\ref{th:main},  these results provide
a way of computing Hausdorff measures of curves as well as a geometric
interpretation of Hausdorff and infinitesimal measures. 
\erem

\bc\label{cor:aux}
Let $\g:[a,b]\rightarrow M$ be a \kdiff\ curve.  Then, for every
measurable set $A\subset [a,b]$, 
$$
\ss^k(\g(A))= {\cal H}^k(\g(A)) ~~\mathrm{ and
}~~\hh^k(\g(A))\leq\Len_k(\g(A)). 
 $$
 If moreover $\g$ is injective then
 $$
 \hh^k(\g(A))=\Len_k(\g(A)).
 $$
 \ec

\brem
Let us consider the case where $\g$ is an injective \kdiff\
curve. Recalling the 
definition of $\Len_k$, we have 
$$
\hh^{k}(\cc)=\int_a^b\meas_t^{k}(\g)\, dt,
$$
that is, we have an integral formula for the $k$-dimensional Hausdorff
measure. 
  Moreover,  Theorem~\ref{th:main} implies that the Hausdorff
  dimension $k_\hh$ of $\cc$ coincides with the length dimension of
  $\cc$. If in addition $\hh^{k_\hh}(\cc)$ (or $\Len_{k_\hh}(\cc)$) is
  finite, then Corollary~\ref{cor:aux} implies that
  $\hh^{k_\hh}\lfloor_\cc$ is absolutely continuous with respect to
  the push-forward   measure\footnote{Given a Borel set 
$E\subset M$ the push-forward measure $\g^*{\cal L}^1$ is defined by
$$\g_*{\cal L}^1(E)={\cal L}^1(\g^{-1}(E\cap \cc)).
$$}
$\g_*{\cal L}^1$ and that its  Radon--Nikodym derivative is
$\meas_t^k(\g)$.
\erem
\medskip

\noindent{\bf Proof of Theorem~\ref{th:main}.} 
Clearly it suffices to prove the first statement  of the theorem.
Consider the (possibly empty) open subset of $[a,b]$
\begin{equation}\label{eq:I}
I = \{ t \in [a,b] \ : \ \meas_t^k(\g)\neq 0
\},
\end{equation}
and its  complementary $I^c=[a,b] \setminus I $. The set $I$ is the
union of a disjointed countable family of  open subintervals $I_i$ of
$[a,b]$. Note that, since $\mathrm{meas}_t^k (\gamma) = 0$ for all $t\in
I^c$, one has
$$
\Len_k (\cc) = \int_I \mathrm{meas}_t^k (\gamma) dt = \sum_i
\int_{I_i} \mathrm{meas}_t^k (\gamma) dt.
$$
%
By Remark~\ref{rk:gen}, we have the equality
\begin{equation}
\label{eq:hkIj}
\hh^k (\gamma(I_i)) = \int_{I_i} \mathrm{meas}_t^k (\gamma) dt,
\quad \forall i.
\end{equation}
Since $\g$ is injective, we have
$$
\hh^k(\cc) \geq \sum_i \hh^k (\gamma(I_i)),
$$
whence we obtain
$\hh^k(\cc) \geq \Len_k (\cc)$.
\medskip

The next step is to prove the converse inequality. Let $\delta>0$.
Since the function $t \mapsto \meas_t^k(\g)^{1/k}$ is uniformly continuous
on $[a,b]$, there exists $\eta>0$ such that, if $t,t' \in [a,b]$ and
$|t-t'|<\eta$, then $|\meas_t^k(\g)^{1/k}-\meas_{t'}^k(\g)^{1/k}|<\delta$. In the covering
 $I=\bigcup_i I_i$, only a finite number
$N_\delta$ of subintervals $I_i$ may have a Lebesgue measure
greater than $\eta$. Up to reordering, we assume $\ll^1 (I_i) <
\eta$ if $i>N_\delta$. Set $J= I^c \cup \bigcup_{i>N_\delta} I_i$.
Since the restriction  of $\meas^k_t(\g)$ to $I^c$ is identically
zero, there holds $\meas_t^k(\g)^{1/k} < \delta$ for every $t\in J$.

The $k$-dimensional Hausdorff measure of $\cc$ satisfies
\begin{equation}
\label{eq:sumIjJ}
\hh^k(\cc) \leq \sum_{i\leq N_\delta} \hh^k(\gamma(I_i)) +
 \hh^k(\gamma(J))
= \sum_{i\leq
N_\delta} \int_{I_i} \mathrm{meas}_t^k (\gamma) dt +
\hh^k(\gamma(J)),
\end{equation}
in view of~(\ref{eq:hkIj}).

It remains to compute $\hh^k(\gamma(J))$. Being the complementary of
$\bigcup_{i\leq N_\delta} I_i$ in $[a,b]$, $J$ is the disjointed
union of $N_\delta +1$ closed subintervals $J_i=[a_i,b_i]$ of
$[a,b]$. For each one of these intervals we will proceed as in the proof
of Proposition~\ref{p:main}.

Let $\eps>0$ and $i \in \{1, \dots , N_\delta +1\}$. We denote by
$N'$ the smallest integer such that $b_i-a_i \leq N'
(\frac{\eps}{2\delta})^k$ and define $t_0, \dots, t_{N'}$ by
$$
t_j = a_i + j \big(\frac{\eps}{2\delta}\big)^k \quad \hbox{for } j=
0, \dots, N'-1, \qquad t_{N'} = b_i.
$$
We then set $S_j = \gamma([t_j,t_{j-1}])$. Applying
Lemma~\ref{le:holder}, we get, for any $t,t' \in [t_j,t_{j-1}]$,
$$
d \big( \gamma(t),\gamma (t') \big) = |t-t'|^{1/k} (\meas_t^k(\g)^{1/k}+  \epsilon_t(t-t')).
$$
Note that $\meas_t^k(\g)^{1/k} < \delta$ since $t\in J$. Note also that, if $\eps$
is small enough, then $\epsilon_t(|t-t'|)$ is
smaller than $\delta$. Therefore $d \big( \gamma(t),\gamma (t')
\big) < 2\delta |t-t'|^{1/k}\leq \eps$ and $\diam S_j \leq \eps$.
As a consequence
$$
\hh^k_\eps (\gamma(J_i)) \leq N' \eps^k \leq (2\delta)^k (b_i-a_i) +
\eps^k,
$$
and $\hh^k (\gamma(J_i)) \leq (2\delta)^k (b_i-a_i)$. It follows that
$$
\hh^k(\gamma(J)) \leq  \sum_{i\leq N_\delta +1} (2\delta)^k
(b_i-a_i) \leq (2\delta)^k (b-a).
$$
Finally, formula~(\ref{eq:sumIjJ}) yields
$$
 \hh^k(\cc) \leq \sum_{i\leq
N_\delta} \int_{I_i} \mathrm{meas}_t^k (\gamma) dt + (2\delta)^k
(b-a).
$$
Letting $\delta \to 0$, we get $\hh^k(\cc) \leq \int_I
\mathrm{meas}_t^k (\gamma) dt = \Len_k (\cc)$, and thus $\hh^k(\cc)
= \Len_k (\cc)$.
Similarly we can show that $\ss^k(\cc)$ and the limit of
$\eps^k \sigma_{\mathrm{int}} (\cc,
\eps)$ are equal to
$\Len_k (\cc)$.
\hfill$\blacksquare$

\medskip

\noindent{\bf Proof of Corollary~\ref{cor:aux}.} When $\g$ is
injective, the conclusions follow from Theorem~\ref{th:main} and from
the regularity of ${\cal L}^1$ and  $\hh^k$ measures (see
Remark~\ref{rk:gen}). 

Assume now that  $\g$ is not injective. We slightly modify
 the second part of the proof of Theorem~\ref{th:main}  replacing the equality in  \r{eq:sumIjJ}
 by
 $$
 \hh^k(\cc) \leq \sum_{i\leq N_\delta} \hh^k(\gamma(I_i)) +
 \hh^k(\gamma(J))
\leq \sum_{i\leq
N_\delta} \int_{I_i} \mathrm{meas}_t^k (\gamma) dt +
\hh^k(\gamma(J)),
 $$
 which is a consequence of Corollary~\ref{cor:leq}. This shows that $\hh^k(\cc)\leq \Len_k(\cc)$ and therefore  $\hh^k(\g(A))\leq\Len_k(\g(A))$.
 Moreover, we have $\hh^k(\g(I^c))=0$, where $I^c$ is the complementary of the set $I$ defined in \r{eq:I}, which in turn implies $\ss^k(\g(I^c))=0$. Thus
 $$
 \ss^k(\cc)=\ss^k(\g(I))=\hh^k(\g(I))=\hh^k(\cc),
 $$
 where the second equality results from Corollary~\ref{cor:leq}.
\hfill$\blacksquare$

\subsection{Generalization to non \kdiff\  curves}
\label{sec:tg}
In this section we present some possible generalizations
of the preceding results (in particular
Theorem~\ref{th:main})  to non \kdiff\ curves.

Consider first the case of  a  continuous curve  $\gamma: [a,b] \to
M$, $\cc=\gamma([a,b])$. For $k\geq 1$, we define
 $I^k$ to be the set of points $t\in[a,b]$ such that $\meas_t^k(\g)$
 is not continuous at $t$ (that is,  such that $\g$ is not \kdiff\ at
 $t$).  
 A standard argument of  measure theory allows to show the following fact.
Assume that  $[a,b] \setminus I^k$ is an open subset of $[a,b]$
 of full ${\cal L}^1$ measure and that $\hh^{k}(\gamma(I^k))=0$. Then
 the conclusions of Corollary~\ref{cor:aux} still hold. Moreover, if $\g$ is injective, then 
 the equalities between $\hh^{k}(\cc)$,
 $\ss^{k}(\cc)$, and $\Len_{k} (\cc)$ as in Theorem~\ref{th:main} hold true. The
 result on the complexity is not valid anymore, since 
$\lim_{\eps\rightarrow 0}\eps^k\sigma_{\mathrm{int}}(\cc,\eps)$ is not
a measure.
  A curve $\cc$ satisfying the properties above actually appears
as a particular case of $({\cal
  H}^k,1)$-rectifiable set, which will be studied in the next section.

It is however worth to mention a consequence of  the result claimed above (and of
Proposition~\ref{ex:uno}) in the context of Carnot--Carath\'eodory
spaces. Let $(M,\bD,g)$ be a sub-Riemannian manifold,  $\gamma: [a,b]
\to M$ be an absolutely continuous injective  curve, and 
$\cc=\gamma([a,b])$. Let $m_\cc\geq 1$ be the smallest integer
such that $\dot\gamma(t)\in \bD^{m_\cc}(\g(t))$ almost everywhere. We denote
 by $I_\cc$ the set of points $t\in[a,b]$ such that either $\gamma$ is
 not $\con^1$ at $t$ or $\gamma(t)$ is 
$\cc$-singular.
\begin{corollary}\label{co:ac}
Assume that $[a,b] \setminus I_\cc$ is an open subset of $[a,b]$
 of full ${\cal L}^1$ measure and $\hh^{m_\cc}(\gamma(I_\cc))=0$.
Then, for any $k\geq 1$,
$$
 \hh^k(\cc) = \ss^k(\cc)=  \Len_k (\cc), \qquad \hbox{ and } \qquad
 \dim_\hh \cc = m_\cc.
$$
\end{corollary}
When the sub-Riemannian manifold is equiregular, it is already
known~\cite[p.\ 104]{gro96} that the
Hausdorff dimension of a one-dimensional submanifold $C$ is the
smallest integer $k$ such that $T_q C \subset\bD^k ( q )$ for every $q
\in C$. Corollary~\ref{co:ac}  generalizes  this fact.\medskip

Any injective \kdiff\ curve being bi-H\"{o}lder  (see
Lemma~\ref{le:holder}), it is also natural to generalize our results
to such curves.  Consider then a bi-H\"{o}lder curve of exponent $1/k$,
i.e. a curve $\gamma: 
[a,b] \to M$, i.e.,  there
exist positive constants $\delta_-$ and $\delta_+$  such that, for
every $t,t+s \in [a,b]$, 
$$
\delta_-|s|^{1/k} \leq d \big( \gamma(t),\gamma (t+s)
\big) \leq
\delta_+|s|^{1/k}.
$$
For such a curve $\cc=\gamma([a,b])$, the $k$-dimensional length does
not always exists but 
Lemma~\ref{le:bounds} gives estimates of  $\hh^{k}(\cc)$ and
 $\ss^{k}(\cc)$ in function of $T=b-a$, and the following  result for
 upper and lower density:
\begin{equation}
\label{eq:density_bihol}
 \left(\frac{\delta_-}{\delta_+}\right)^{k} \leq \liminf_{
 r \rightarrow 0^+
}\frac{{\cal H}^k(\cc\cap B(q,r))}{2r^k}\leq \limsup_{
 r \rightarrow 0^+
}\frac{{\cal H}^k(\cc\cap B(q,r))}{2r^k}\leq 
\left(\frac{\delta_+}{\delta_-}\right)^{k}. 
\end{equation}
Let us remark  that there is no hope to obtain a density result such
as~\eqref{eq:den} for bi-H\"{o}lder curves. Indeed Assouad
proved in~\cite{assouad} that, for any $k<n$,  there exist
bi-H\"{o}lder curves of exponent $1/k$ from $(-1,1)$ to $\R^n$ (both
endowed with a 
Euclidean metric). When $k>1$ the density of these curves cannot be
constant, for otherwise Lemma~\ref{le:riem_den} would yield a
contradiction. This strongly hints that there is not a
Rademacher's-type result  
 in this context, that is,  being bi-H\"{o}lder of exponent $1/k$ does
 not imply being 
\kdiff\ almost everywhere.

In what follows, in particular in the definition of $({\cal
  H}^k,1)$-rectifiability, we will work with \kdiff\ curves and not
with bi-H\"{o}lder curves. The drawback is that our definitions will not be 
invariant under bi-Lipschitz equivalence of metric spaces since the
\kdiff\ property is not invariant under such equivalence, contrarily
to the bi-H\"{o}lder property. However we think that 
rectifiable sets should be defined as sets which admit almost
everywhere a metric derivative.  As noticed above, the use of
bi-H\"{o}lder curves would not guarantee such a property.

\section{$({\cal H}^k,1)$-rectifiable sets and a density result}
\label{sec:rect}
In this section we use \kdiff\ curves to define a new class of
$(\hh^k,1)$-rectifiable subsets of metric spaces. 

Consider the Euclidean space $\R^n$. Recall that, given  a positive
measure $\mu$ on Borel subsets of $\R^n$, a subset $S\subset \R^n$ is
\emph{$(\mu,k)$-rectifiable} if there exists 
a countable family of Lipschitz functions $\g_i:V_i\rightarrow \R^n$,
$i\in\N$, 
where $V_i$ is a bounded  subset of $\R^{k}$, such that
$\mu(S\setminus\cup_{i\in\N}\g_i(V_i))=0$ (see
\cite[3.2.14]{federer}). Considering $\mu=\hh^{k'}$ on $\R^n$ (with
the Euclidean structure), one has that if $k'>k$ then there are no
$(\hh^{k'},k)$-rectifiable sets of positive $\hh^{k'}$ measure. This
follows from the requirement of $\g_i$ being  Lipschitz. 
On the other hand, if one consider images under H\"older continuous
functions $\g_i$ then the case $k'>k$ becomes of interest. This
suggests the next definition. 

 Consider a metric space $(M,d)$.

\bdeff\label{def:rect}
A subset $S\subset M$ is $({\cal H}^k,1)$\emph{-rectifiable} if there exists
 a countable family of \kdiff\   curves $\g_i:I_i\rightarrow M$,
 $i\in\N$, $I_i$ closed interval in $\R$
such that
$$
{\cal H}^k(S\setminus \cup_{i\in\N}\g_i(I_i))=0.
$$
\edeff
\brem
If $M$ is a manifold and $d$ is the distance associated with a Riemannian structure on $M$, the class of $({\cal H}^k,1)$-rectifiable sets with positive and finite $\hh^k$ measure is empty unless $k=1$ (see Proposition~\ref{prop:eucl}).  Since m-$\con^1_1$ curves   are Lipschitz, 
 in this case Definition~\ref{def:rect} coincides with the usual definition of $({\cal H}^1,1)$-rectifiable sets. Conversely, when $(M,d)$ is the Carnot--Carath\'eodory space associated with a genuine sub-Riemannian manifold,  there exist $({\cal H}^k, 1)$-rectifiable sets (of positive and finite $\hh^k$ measure) for some integers $k>1$ (see Section~\ref{sec:examples}). \erem
When a subset is $\hh^k$-measurable and has finite $\hh^k$ measure, being $(\hh^k,1)$-rectifiable implies boundedness for the lower and upper densities of the measure $\hh^k\lfloor_S$.
\bt\label{th:fico}
Assume $S\subset M$ is a $({\cal H}^k,1)$-rectifiable and $\hh^k$-measurable set such that $\hh^k(S)<+\infty$.
Then  for ${\cal H}^k$-almost every $q\in S$
\begin{equation}\label{eq:density}
2\leq \liminf_{
 r \rightarrow 0^+
}\frac{{\cal H}^k(S\cap B(q,r))}{r^k}\leq \limsup_{
 r \rightarrow 0^+
}\frac{{\cal H}^k(S\cap B(q,r))}{r^k}\leq 2^{k}.
\end{equation}
\et
%

Recall that in \cite[Co. 5.5]{preiss} it was proved that, in the
Euclidean case, there exists a constant $c>0$ such that if a
$\mu$-measurable subset $E\subset \R^n$  with finite $\mu$ measure
satisfies
\begin{equation}\label{eq:preiss}
0< \limsup_{
 r \rightarrow 0^+
}\frac{\mu(E\cap B(x,r))}{r^k}\leq c \liminf_{
 r \rightarrow 0^+
}\frac{\mu(E\cap B(x,r))}{r^k}<+\infty,
\end{equation}
 for $\mu$-almost every $x\in E$,
then $E$ is $(\mu,k)$-rectifiable. This result provides a
characterization of rectifiable sets as the converse is also true (see
\cite[Th. 3.2.19]{federer}).
 Theorem~\ref{th:fico} implies  that if $S\subset M$ is
 $(\hh^k,1)$-rectifiable in the sense of Definition~\ref{def:rect}
 then
 \begin{equation}\label{eq:oq}
 0< \limsup_{
 r \rightarrow 0^+
}\frac{{\cal H}^k(S\cap B(q,r))}{r^k}\leq 2^{k-1}\liminf_{
 r \rightarrow 0^+
}\frac{{\cal H}^k(S\cap B(q,r))}{r^k}<+\infty,
 \end{equation}
 for $\hh^k$-almost every $q\in S$. The last estimate is, mutatis mutandis,  the assumption \r{eq:preiss} in the result by Preiss.
An open question is whether the same conclusion of \cite[Cor.\ 5.5]{preiss} holds  with our notion of $(\hh^k,1)$-rectifiable sets. Namely,  is condition \r{eq:oq} for $\hh^k$-almost every $q\in S\subset M$ sufficient to show that a $\hh^k$-measurable set $S$ of finite $\hh^k$ measure is $(\hh^k,1)$-rectifiable in the sense of Definition~\ref{def:rect}?

\smallskip

\noindent{\bf Proof of Theorem~\ref{th:fico}.}
By assumption, there exists a countable family of \kdiff\ curves $\g_i:I_i\to M$ such that $I_i$ is a closed interval and $\hh^k(S\setminus \cup_i\g_i(I_i))=0$.  Since by Corollary~\ref{cor:aux}, for every $i$, $\hh^k(\g_i(\{t\mid \meas^k_t(\g_i)=0\}))=0$, we may assume $\meas^k_t(\g_i)\neq 0$ for every $t\in I_i$ and then, by a reparameterization, $\meas^k_t(\g_i)\equiv 1$. %
%
This implies that every $\g_i$ is locally injective. Hence,
without loss of generality, we may  assume that every $\g_i$ is injective and moreover that the sets $\g_i(I_i)$ are pairwise disjoint.

\smallskip

Since $\hh^k(S)<+\infty$, to prove the upper bound
$$
 \limsup_{
 r \rightarrow 0^+
}\frac{{\cal H}^k(S\cap B(q,r))}{r^k}\leq 2^{k},
$$
for $\hh^k$-almost every $q\in M$, it suffices to use \cite[2.10.19 (5)]{federer}.

Let us show the lower bound in \r{eq:density}, namely, that
\begin{equation}\label{eq:liminf}
\liminf_{
 r \rightarrow 0^+
}\frac{{\cal H}^k(S\cap B(q,r))}{r^k}\geq 2,
\end{equation}
for $\hh^k$-almost every $q\in S$.
Let $\tilde I_i=\g_i^{-1}(\g_i(I_i)\cap S)$. Then $\cup_{i\in N}\g_i(\tilde I_i)\subset S$ and $\hh^k(S\setminus\cup_{i\in N}\g_i(\Tt I_i))=0$.
 We may assume $\hh^k(\g_i(\Tt I_i))>0$ for each $i$. Then, by Corollary~\ref{cor:aux}, since $\meas^k_t(\g_i)\equiv 1$,
${\cal L}^1(\Tt I_i)=\hh^k(\g_i(\Tt I_i))> 0$. Therefore almost every $t\in\Tt I_i$ is a density point for the Lebesgue measure on $\Tt I_i$, i.e.,
$$
\lim_{r\to 0}\frac{{\cal L}^1(\Tt I_i\cap B(t,r))}{2 r}=1,
$$
where $B(t,r)=(t-r,t+r)$.
Hence, for $\hh^k$-almost every $q\in S$ there exist a unique $i$ and a unique $t\in \Tt I_i$ such that $q=\g_i(t)$ and $t$ is a density point for ${\cal L}^1\lfloor_{\Tt I_i}$.
Since $\g_i(\Tt I_i)\subset S$, we deduce
$$
\frac{\hh^k(S\cap B(q,r))}{r^k}\geq \frac{\hh^k(\g_i(\Tt I_i)\cap B(q,r))}{r^k}=\frac{{\cal L}^1(\Tt I_i\cap \g_i^{-1}(B(q,r)))}{r^k},
$$
the last equality following by Corollary~\ref{cor:aux}.
Now, for any $\delta >0$, from
Lemma~\ref{le:holder}, for $|t-s|\leq \frac{r^k}{(1+\delta)^k}$ we have
$$
d \big( \gamma(t),\gamma (s)
\big) \leq |t-s|^{1/k}(1+\delta)\leq r.
$$
This implies $B(t,r^k/(1+\delta)^k)\subset \g_i^{-1}(B(q,r))$. Therefore
$$
\frac{{\cal L}^1(\Tt I_i\cap \g_i^{-1}(B(q,r)))}{r^k}\geq \frac{{\cal L}^1(\Tt I_i\cap B(t,r^k/(1+\delta)^k)))}{r^k}.
$$
The right-hand side of the inequality above tends to $2/(1+\delta)^k$, as $r$ goes to $0$, since $t$ is a density point for ${\cal L}^1\lfloor_{\Tt I_i}$. Letting $\delta$ go to $0$, we conclude
$$
\liminf_{r\to 0}\frac{{\cal L}^1(\Tt I_i\cap \g_i^{-1}(B(q,r)))}{r^k}\geq 2,
$$
which shows \r{eq:liminf}.
\hfill$\blacksquare$

\bibliographystyle{abbrv}
\bibliography{biblio_mdiff}

\begin{thebibliography}{10}

\bibitem{balu}
A.~Agrachev, D.~Barilari, and U.~Boscain.
\newblock On the {H}ausdorff volume in sub-{R}iemannian geometry.
\newblock preprint arXiv:1005.0540v3.

\bibitem{ambrosio-rectifiable}
L.~Ambrosio and B.~Kirchheim.
\newblock Rectifiable sets in metric and {B}anach spaces.
\newblock {\em Math. Ann.}, 318(3):527--555, 2000.

\bibitem{ambrosio}
L.~Ambrosio and P.~Tilli.
\newblock {\em Topics on analysis in metric spaces}, volume~25 of {\em Oxford
  Lecture Series in Mathematics and its Applications}.
\newblock Oxford University Press, Oxford, 2004.

\bibitem{assouad}
P.~Assouad.
\newblock Plongements lipschitziens dans {${\bf R}^{n}$}.
\newblock {\em Bull. Soc. Math. France}, 111(4):429--448, 1983.

\bibitem{bellaiche}
A.~Bella{\"{\i}}che.
\newblock The tangent space in sub-{R}iemannian geometry.
\newblock In {\em Sub-{R}iemannian geometry}, volume 144 of {\em Progr. Math.},
  pages 1--78. Birkh\"auser, Basel, 1996.

\bibitem{bes1}
A.~S. Besicovitch.
\newblock On the fundamental geometrical properties of linearly measurable
  plane sets of points.
\newblock {\em Math. Ann.}, 98(1):422--464, 1928.

\bibitem{lipeq}
U.~Boscain, G.~Charlot, R.~Ghezzi, and M.~Sigalotti.
\newblock {L}ipschitz classification of two-dimensional almost-{R}iemannian
  distances on compact oriented surfaces.
\newblock preprint arXiv:1003.4842, to appear on {\it Journal of Geometric
  Analysis}, 2011.

\bibitem{j-falbel}
E.~Falbel and F.~Jean.
\newblock Measures of transverse paths in sub-{R}iemannian geometry.
\newblock {\em J. Anal. Math.}, 91:231--246, 2003.

\bibitem{fed}
H.~Federer.
\newblock The {$(\varphi,k)$} rectifiable subsets of {$n$}-space.
\newblock {\em Trans. Amer. Soc.}, 62:114--192, 1947.

\bibitem{federer}
H.~Federer.
\newblock {\em Geometric measure theory}.
\newblock Die Grundlehren der mathematischen Wissenschaften, Band 153.
  Springer-Verlag New York Inc., New York, 1969.

\bibitem{rect-heis}
B.~Franchi, R.~Serapioni, and F.~Serra~Cassano.
\newblock Rectifiability and perimeter in the {H}eisenberg group.
\newblock {\em Math. Ann.}, 321(3):479--531, 2001.

\bibitem{gz-nuovo}
J.-P. Gauthier, B.~Jakubczyk, and V.~Zakalyukin.
\newblock Motion planning and fastly oscillating controls.
\newblock {\em SIAM J. Control Optim.}, 48(5):3433--3448, 2009/10.

\bibitem{gau06}
J.-P. Gauthier and V.~Zakalyukin.
\newblock On the one-step-bracket-generating motion planning problem.
\newblock {\em J. Dyn. Control Syst.}, 11(2):215--235, 2005.

\bibitem{gro96}
M.~Gromov.
\newblock Carnot-{C}arath\'eodory spaces seen from within.
\newblock In {\em Sub-{R}iemannian geometry}, volume 144 of {\em Progr. Math.},
  pages 79--323. Birkh\"auser, Basel, 1996.

\bibitem{hardy}
G.~H. Hardy.
\newblock Weierstrass's non-differentiable function.
\newblock {\em Trans. Amer. Math. Soc.}, 17(3):301--325, 1916.

\bibitem{jea00}
F.~Jean.
\newblock {\em Paths in Sub-Riemannian Geometry}.
\newblock Springer (A. Isidori, F. Lamnabhi-Lagarrigue and W. Respondek Eds.),
  2000.

\bibitem{cocv_fj}
F.~Jean.
\newblock Entropy and complexity of a path in sub-{R}iemannian geometry.
\newblock {\em ESAIM Control Optim. Calc. Var.}, 9:485--508 (electronic), 2003.

\bibitem{kircheim}
B.~Kirchheim.
\newblock Rectifiable metric spaces: local structure and regularity of the
  {H}ausdorff measure.
\newblock {\em Proc. Amer. Math. Soc.}, 121(1):113--123, 1994.

\bibitem{magnani}
V.~Magnani.
\newblock Characteristic points, rectifiability and perimeter measure on
  stratified groups.
\newblock {\em J. Eur. Math. Soc.}, 8(4):585--609, 2006.

\bibitem{marstrand}
J.~M. Marstrand.
\newblock The {$(\varphi ,\,s)$} regular subsets of {$n$}-space.
\newblock {\em Trans. Amer. Math. Soc.}, 113:369--392, 1964.

\bibitem{mattila}
P.~Mattila.
\newblock {\em Geometry of sets and measures in Euclidean spaces: fractals and
  rectifiability}.
\newblock Cambridge Studies in Advanced Mathematics. Cambridge University
  Press, 1996.

\bibitem{rect-heis2}
P.~Mattila, R.~Serapioni, and F.~Serra~Cassano.
\newblock Characterizations of intrinsic rectifiability in {H}eisenberg groups.
\newblock {\em Ann. Sc. Norm. Super. Pisa Cl. Sci. (5)}, 9(4):687--723, 2010.

\bibitem{moore}
E.~F. Moore.
\newblock Density ratios and {$(\phi,1)$} rectifiability in {$n$}-space.
\newblock {\em Trans. Amer. Math. Soc.}, 69:324--334, 1950.

\bibitem{preiss}
D.~Preiss.
\newblock Geometry of measures in {${\bf R}^n$}: distribution, rectifiability,
  and densities.
\newblock {\em Ann. of Math. (2)}, 125(3):537--643, 1987.

\bibitem{weier}
K.~Weierstrass.
\newblock {\em On continuous functions of a real argument that do not have a
  well-defined differential quotient}.
\newblock G.A. Edgar, Classics on Fractals. Addison-Wesley Publishing Company,
  1993.

\end{thebibliography}

\end{document}